\documentclass[10pt]{article}
\usepackage{latexsym}
\usepackage[dvips]{color}
\usepackage{float}
\usepackage{afterpage} 
\usepackage [pdftex]{graphicx}
\usepackage{epsfig}
\usepackage{amsmath}
\usepackage{amsfonts}
\usepackage{amssymb}
\usepackage{setspace}
\usepackage{subfigure}

\title{Using Monoidal Categories in the Transformational Study of Musical Time-Spans and Rhythms}
\author{
        Alexandre Popoff \\
                al.popoff@free.fr\\
        France}

\begin{document}

\maketitle

\section{Introduction}

The field of musical transformational theory, initiated by the work of David Lewin \cite{lewin}, relies on the use of a group structure, in which group elements can be seen as transformations (or as Lewin called them, "generalized intervals") between musical objects. In neo-Riemannian theory, the typical set of elements is constituted by the major and minor chords, and the most famous corresponding group of transformations is the group isomorphic to the dihedral group $D_{24}$ of 24 elements. This group can act on the set of major and minor chords through the well-known L, R and P operations or through the transpositions and inversions operators \cite{cohn1,cohn2,cohn3,capuzzo}, though many other actions are also possible (see for example the Schritt-Wechsel group) \cite{douthett}. Other groups \cite{straus,hook1,hook2,peck1,peck2,fiore1,fiore2,fiore3} have been used for the transformation of chords, or more generally pitch-based musical objects. In a general setting, we have shown \cite{popoff} that generalized neo-Riemannian groups of musical transformations can be often be built as extensions.

Comparatively less work has been done regarding the transformational theory of musical time-spans. We recall the definition of a time-span:

\begin{description}
\item[Definition]{\textit{A time-span is a half-open interval of $\mathbb{R}$ of the form $[t,t+\Delta[$, with $t \in \mathbb{R}$ and $\Delta \in \mathbb{R}_+^*$ ($\mathbb{R}_+^*$ being the set of strictly positive real numbers). The value $t$ is called the onset of the time-span, whereas $\Delta$ is its duration. A time-span is equivalently referred to by the pair $(t,\Delta)$.}}
\end{description}

Again, Lewin (\cite{lewin}, pp. 60-81) was the first to provide a group structure and a group action for time-spans. In our previous work (see Section 6 of \cite{popoff}), we have shown that this group structure is also covered by the theory of extensions as applied to musical transformations. More precisely, the group of transformations $G$ is built as an extension of a group $O$ by a group $D$, where $O$ is the group structure of onsets and $D$ is the group structure of durations (making $O$ a normal subgroup of $G$). Some authors have considered the transformational theory of time-series (\cite{morris}), while others also considered the problem of rhythm though not necessarily from a transformational point of view (\cite{hook4,agmon,amiot}).

\begin{figure}
\centering
\includegraphics[scale=0.5]{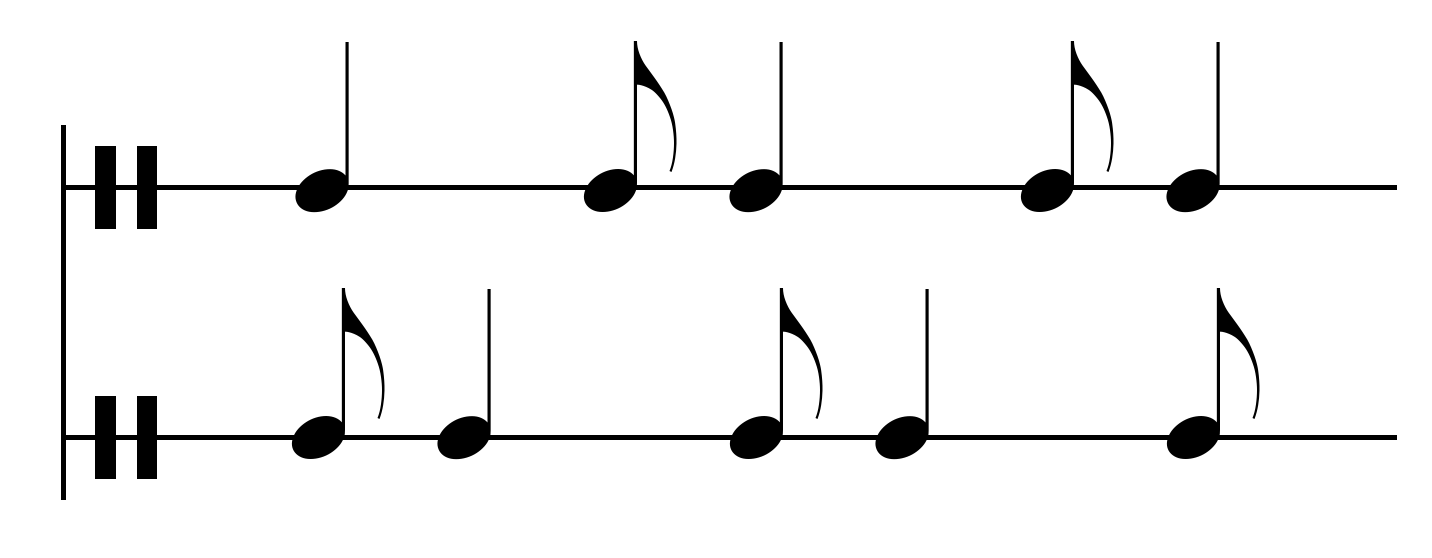}
\caption{A two-parts alternating rhythm.}
\label{fig:Rhythm}
\end{figure}

Our work has shown that the framework of group extensions allows for generalizations of the original transformation group of Lewin. Consider for example the rhythm presented in Figure \ref{fig:Rhythm}.
There is an obvious symmetry between parts, which will not be apparent if we study the transformation of the time-spans in each timeline separately.

In the one-dimensional case, the space of onsets is $\mathbb{R}$ equipped with an additive group structure $(\mathbb{R},+)$, whereas the space of durations is $\mathbb{R}_+^*$ equipped with a multiplicative group structure $(\mathbb{R}_+^*,\times)$. The corresponding group extension is therefore the semidirect product $(\mathbb{R},+) \rtimes (\mathbb{R}_+^*,\times)$, which is a subgroup of the affine group in one dimension. Generalizing to two dimensions, the space of onsets is the direct product $\mathbb{R}^2$ which can be equipped with an additive group structure $(\mathbb{R},+) \times (\mathbb{R},+)$. Intuitively, we could then choose a subgroup of $GL(2,\mathbb{R})$ for the group of durations, so that the final group transformations would be a subgroup of the affine group in two dimensions. However, not every matrix of $GL(2,\mathbb{R})$ can be used; indeed, how can we define an element of $d \in D$ of $GL(2,\mathbb{R})$ to be a "generalized duration" and how can we avoid "negative durations" ?

In the one-dimensional case, the right action of an element $(u,\delta)$ of $(\mathbb{R},+) \rtimes (\mathbb{R}_+^*,\times)$ on a time-span $(t,\Delta)$ is given by

$$ (t,\Delta) \cdot (u,\delta) = (t+\Delta u, \delta \Delta)$$
In particular, the right action of the translation operator $(1,1)$ results in 

$$ (t,\Delta) \cdot (1,1) = (t+\Delta, \Delta)$$
i.e the time-span is translated by an amount of time equal to its duration.

If we apply the same reasoning to the two-dimensional case, the right action of the translation operator $(\left|\begin{matrix}1 \\ 1 \end{matrix}\right|, \left|\begin{matrix} 1 & 0 \\ 0 & 1 \end{matrix}\right|)$ on a time-span $(\left|\begin{matrix}t_1 \\ t_2 \end{matrix}\right|, \left|\begin{matrix} a & b \\ c & d \end{matrix}\right|)$ will result in 

$$ (\left|\begin{matrix}t_1 \\ t_2 \end{matrix}\right|, \left|\begin{matrix} a & b \\ c & d \end{matrix}\right|) \cdot (\left|\begin{matrix}1 \\ 1 \end{matrix}\right|, \left|\begin{matrix} 1 & 0 \\ 0 & 1 \end{matrix}\right|) \cdot = (\left|\begin{matrix}t_1 + (a+b) \\ t_2 + (c+d) \end{matrix}\right|, \left|\begin{matrix} a & b \\ c & d \end{matrix}\right|) $$

Thus we see that for an element of $\left|\begin{matrix} a & b \\ c & d \end{matrix}\right|$ of $GL(2,\mathbb{R})$ the duration on the first time-line is $(a+b)$, and the duration on the second time-line is  $(c+d)$. In order to have strictly positive duration we thus need find a subgroup of matrices of $GL(2,\mathbb{R})$ such that $a+b > 0, c+d>0$.

This reason prompted us to choose in \cite{popoff} a subgroup $U$ of $GL(2,\mathbb{R})$ consisting of all 2$\times$2 matrices $M$ of the form $\left|\begin{matrix}\alpha & 0 \\ 0 & \beta \end{matrix}\right|$ or $\left|\begin{matrix}0 & \alpha \\ \beta & 0 \end{matrix}\right|$, with $\alpha \in \mathbb{R}_+^*$, $\beta \in \mathbb{R}_+^*$. We have shown that the rhythm of Figure \ref{fig:Rhythm} is obtained through the successive right multiplication of the group element $\left(\left|\begin{matrix}1 \\ 1 \end{matrix}\right|,\left|\begin{matrix}0 & \frac{1}{2} \\ 2 & 0 \end{matrix}\right|\right)$ acting on the initial time-span $\left(\left|\begin{matrix}0 \\ 0 \end{matrix}\right|,\left|\begin{matrix}1 & 0 \\ 0 & \frac{1}{2} \end{matrix}\right|\right)$. The matrix being $\left|\begin{matrix}0 & \frac{1}{2} \\ 2 & 0 \end{matrix}\right|$ being equal to $\left|\begin{matrix}\frac{1}{2} & 0 \\ 0 & 2 \end{matrix}\right| \cdot \left|\begin{matrix} 0 & 1 \\ 1 & 0 \end{matrix}\right|$ we thus clearly see symmetry between parts in the reversing matrix $\left|\begin{matrix} 0 & 1 \\ 1 & 0 \end{matrix}\right|$.

Notice however that there exists other possible subgroups of $GL(2,\mathbb{R})$ fullfilling the requirement for positive "generalized durations", for example the subgroup consisting of matrices of the form $\left|\begin{matrix}1 & \alpha \\ 0 & 1+\alpha \end{matrix}\right|$ or $\left|\begin{matrix} \alpha & 1 \\ 1+\alpha & 0 \end{matrix}\right|$ with $\alpha \in \mathbb{R}_+^*$. We have not looked at the maximal subgroup of $GL(2,\mathbb{R})$ of matrices $\left|\begin{matrix} a & b \\ c & d \end{matrix}\right|$ such that $a+b > 0, c+d>0$ and this could potentially be an interesting lead for future musical applications. Notice also that this construction can easily be generalized to $n$ timelines, with $O=(\mathbb{R},+)^n$ the direct product of $n$ copies of $(\mathbb{R},+)$ and $D$ being an appropriate subgroup of $GL(n,\mathbb{R})$.

Our attention was drawn instead by a straightforward result linking the generalized group $(\mathbb{R},+)^2 \rtimes U$ and monoidal categories, which we present in the next section.

\section{Monoidal categories and transformations of $n$ timelines}

We refer the reader to \cite{maclane1} for the definitions regarding monoidal categories. In particular, we will make an extensive use of string diagrams as a convenient way to describe morphisms in a monoidal category, which are reviewed in \cite{selinger}.

To link the group introduced in the previous section with monoidal categories, we first observe the following proposition. Note that we define the wreath product $G \wr S_n$ of a group $G$ by the symmetric group $S_n$ on $n$ elements, as the semidirect product of $n$ copies of $G$ and $S_n$, with the natural action of $S_n$ on the $n$ copies.

\begin{description}
\item[Proposition]{\textit{
Let $O=(\mathbb{R},+)^n$ be the direct product of $n$ copies of $(\mathbb{R},+)$. Let $D$ be the subgroup of $GL(n,\mathbb{R})$ generated by the matrices of the form $\left|\begin{matrix} \alpha_1 & 0  & \cdots \\ 0 & \ddots & 0 \\ 0 & \cdots & \alpha_n \end{matrix}\right| \cdot \sigma$ where $\alpha_i > 0, \forall i$ and $\sigma$ is a permutation matrix of size $n$. Then the semidirect product $O \rtimes D$ is isomorphic to the wreath product $Aff_+(\mathbb{R}) \wr S_n$ where $Aff_+(\mathbb{R})$ is the semidirect product $(\mathbb{R},+) \rtimes (\mathbb{R}_+^*,\times)$ and $S_n$ is the symmetric group on $n$ elements.
}}
\vspace{0.2cm}
\item[Proof]{
By an abuse of notation, we denote by $\sigma$ either a permutation matrix of size $n$ or the corresponding group element in the symmetric group $S_n$.
Let $\left|\begin{matrix} \alpha_1 & 0  & \cdots \\ 0 & \ddots & 0 \\ 0 & \cdots & \alpha_n \end{matrix}\right|$ be a matrix, such that $\alpha_i > 0, \forall i$. It can equivalently be represented by an ordered sequence $(\alpha_1, \cdots, \alpha_n)$. We show that every element of $D$ is of the form $\left|\begin{matrix} \alpha_1 & 0  & \cdots \\ 0 & \ddots & 0 \\ 0 & \cdots & \alpha_n \end{matrix}\right| \cdot \sigma$. Indeed, we have $\sigma \cdot \left|\begin{matrix} \alpha_1 & 0  & \cdots \\ 0 & \ddots & 0 \\ 0 & \cdots & \alpha_n \end{matrix}\right| = \left|\begin{matrix} \alpha'_1 & 0  & \cdots \\ 0 & \ddots & 0 \\ 0 & \cdots & \alpha'_n \end{matrix}\right| \cdot \sigma$ where $(\alpha'_1, \cdots, \alpha'_n)$ results from the permutation $\sigma$ applied to $(\alpha_1, \cdots, \alpha_n)$.

Let $(o_1,d_1 \cdot \sigma_1)$ and $(o_2,d_2 \cdot \sigma_2)$ be two elements of $O \rtimes D$. We denote by $(t_1, \cdot , t_n)$ the components of $o_1$, and by $(u_1, \cdot , u_n)$ the components of $o_2$. Similarly, we denote by $(\alpha_1, \cdots, \alpha_n)$ the components of $d_1$ and by $(\beta_1, \cdots, \beta_n)$ the components of $d_2$. Let $\psi$ be the function which sends any element $(o_1,d_1 \cdot \sigma_1) \in O \rtimes D$ to $((t_1,\alpha_1), \cdots , (t_n,\alpha_n) , \sigma_1) \in Aff_+(\mathbb{R}) \wr S_n$. This function is a bijection, and we need to show this is an homomorphism.

The function $\psi$ sends the identity of $O \rtimes D$ to the identity of $Aff_+(\mathbb{R}) \wr S_n$. Furthermore, we have

$$ (o_1,d_1 \cdot \sigma_1)(o_2,d_2 \cdot \sigma_2) = (o_1+d_1 \sigma_1 o_2,d_1 \cdot \sigma_1 \cdot d_2 \cdot \sigma_2)$$

with

$$o_1+d_1 \sigma_1 o_2 = (t_1+\alpha_1 \sigma_1(u_1), \cdots , t_n+\alpha_n \sigma_1(u_n))$$

(where $\sigma_1(u_i)$ is the image of $u_i$ by the permutation $\sigma_1$) and 

$$d_1 \cdot \sigma_1 \cdot d_2 \cdot \sigma_2 = (\alpha_1 \sigma_1(\beta_1), \cdots , \alpha_n \sigma_1(\beta_n)) \sigma_1 \sigma_2.$$

On the other hand, the multiplication of $\psi((o_1,d_1 \cdot \sigma_1))$ and $\psi((o_2,d_2 \cdot \sigma_2))$ in $Aff_+(\mathbb{R}) \wr S_n$ gives 

\begin{multline*}
((t_1,\alpha_1), \cdots , (t_n,\alpha_n) , \sigma_1) \ast ((u_1,\beta_1), \cdots , (u_n,\beta_n) , \sigma_2)\\ =((t_1,\alpha_1) \cdot \sigma_1((u_1,\beta_1)), \cdots , (t_n,\alpha_n) \cdot \sigma_1((u_n,\beta_n)), \sigma_1 \sigma_2)\\
=((t_1,\alpha_1) \cdot (\sigma_1(u_1),\sigma_1(\beta_1)), \cdots , (t_n,\alpha_n) \cdot (\sigma_1(u_n),\sigma_1(\beta_n)), \sigma_1 \sigma_2)\\
=((t_1+\alpha_1\sigma_1(u_1),\alpha_1\sigma_1(\beta_1)), \cdots , (t_n+\alpha_1\sigma_1(u_n),\alpha_n\sigma_1(\beta_n)), \sigma_1 \sigma_2)
\end{multline*}
Hence $\psi$ is the desired isomorphism.
\begin{flushright}$\square$\end{flushright}
}
\end{description}

We then prove the following result

\begin{description}
\item[Proposition]{\textit{
Let $G$ be a group considered as a category $\mathcal{G}$, i.e a category with one object $A$ and invertible morphisms $Hom(A,A)=G$. Let $\cal{C}$ be the free strict symmetric monoidal category generated by the category $\mathcal{G}$. Then $\forall n \in \mathbb{N}$, $Hom(A^{\otimes n},A^{\otimes n})$ is isomorphic to the wreath product $G \wr S_n$.
}}
\vspace{0.2cm}
\item[Proof]{
The construction of the symmetric monoidal category $\cal{C}$ defines a braiding isomorphism $\gamma: A \otimes A \to A \otimes A$, such that $\gamma = \gamma^{-1}$. Therefore all morphisms of $A^{\otimes n}$ are invertible, i.e $Hom(A^{\otimes n},A^{\otimes n})$ is a group. Notice that for a given $n$, the braiding $\gamma$ allows the definition of elementary morphisms $\sigma_i$ in $Hom(A^{\otimes n},A^{\otimes n})$ with $\sigma_i=(id \otimes \cdots \otimes id \otimes \gamma \otimes id \otimes \cdots \otimes id)$ where $\gamma$ is in the $i-$th position ($1<i<n$). The morphisms $\sigma_i$ generates the possible permutations $\sigma$ of the $n$ tensor copies of $A$. Thus, a morphism of $A^{\otimes n}$ can be put in the form $(g_1\otimes \cdots \otimes g_n) \circ \sigma$ where $g_i  \in G$ and $\sigma$ is generated by the $\sigma_i$ ($1<i<n$). Obviously, $\sigma$ can be mapped to an element of $S_n$. Let $\psi$ be the function which sends an element $(g_1\otimes \cdots \otimes g_n) \circ \sigma$ of $Hom(A^{\otimes n},A^{\otimes n})$ to the element $(g_1, \cdots , g_n , \sigma)$ of $G \wr S_n$. It is easily verified that $\psi$ is the desired isomorphism of groups.
\begin{flushright}$\square$\end{flushright}
}
\end{description}

By combining the above propositions we arrive at the following musical interpretation for the transformational study of time-spans on different timelines. The "generalized interval system" (GIS) approach of Lewin has been proved equivalent to the study of simply transitive group actions on sets \cite{fiore4}.
From a categorical approach, it is known that a simply transitive group action on a set is equivalent to having a representable functor from the group G viewed as a category, i.e a category with only one object $A$ and invertible morphisms, to $\mathbf{Set}$. Recall that a covariant representable functor is a functor naturally isomorphic to the $Hom(A,-)$ functor. From the Yoneda lemma, the choice of a representable functor is equivalent to picking one particular element of the set $S$ as the group identity. This defines a bijection $\chi : G \to S$ and a canonical (left) group action given by
$$g \cdot s = \chi(g \cdot \chi^{-1}(s))$$
where $s$ is an element of $S$ and $g$ is an element of $G$. If contravariant functors are used, then a contravariant representable functor will be naturally isomorphic to the $Hom(-,A)$ functor and, after establishing the bijection $\chi$, the canonical (right) group action will be given by
$$s \cdot g = \chi(\chi^{-1}(s) \cdot g)$$

Applying this to time-spans in particular, if we take the category $\mathcal{G}$ with only one object $A$ such that $Hom(A,A) \cong Aff_+(\mathbb{R})$, then the object $A$ is analogous to a single timeline and the morphisms of $A$ can be viewed as the collection of time-spans on this timeline. By constructing the free symmetric monoidal category $\mathcal{C}$ generated by $\mathcal{G}$, then the objects of $\mathcal{C}$ are all the possible tensor powers $A^{\otimes n}$ of $A$, $n \in \mathbb{N}$. Tensoring is interpreted as juxtaposing timelines, and thus $A^{\otimes n}$ is analogous to having $n$ timelines, the symmetric nature of $\mathcal{C}$ reflecting the possible interaction between timelines when applying transformations. In a configuration of $n$ timelines, the possible time-spans are given by the collection of morphisms of $A^{\otimes n}$ for which we have seen that $Hom(A^{\otimes n},A^{\otimes n})=Aff_+(\mathbb{R}) \wr S_n$.

\begin{figure}
\centering
\subfigure[]{
\includegraphics[scale=0.3]{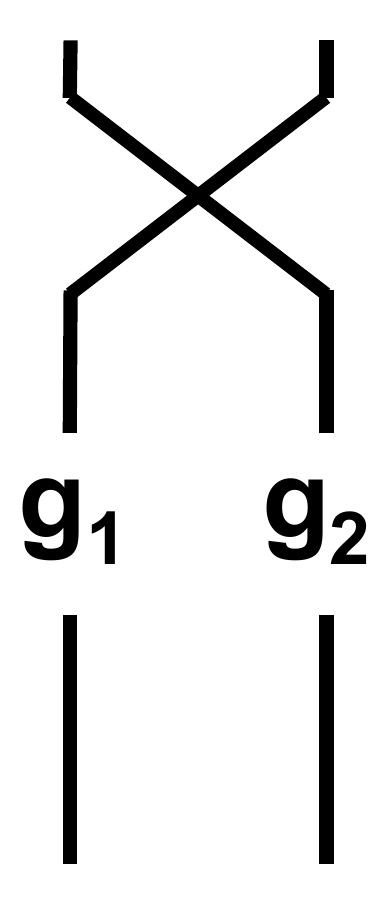}
\label{subfig:generalform}
}
\qquad
\subfigure[]{
\includegraphics[scale=0.3]{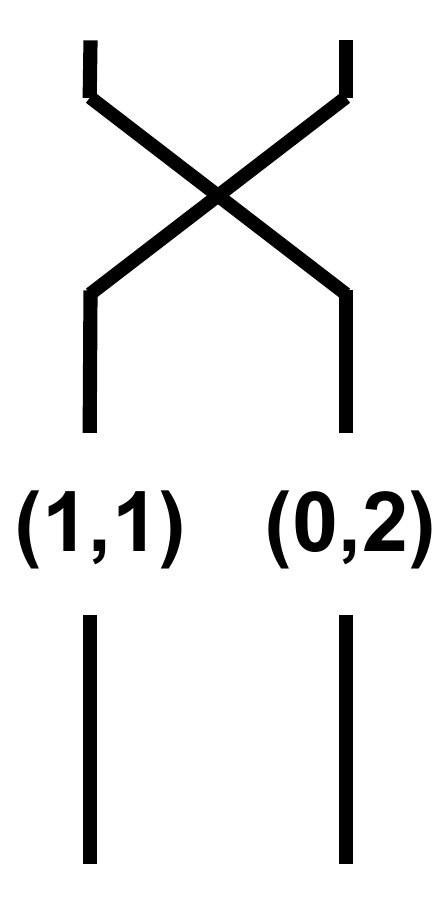}
\qquad
\includegraphics[scale=0.3]{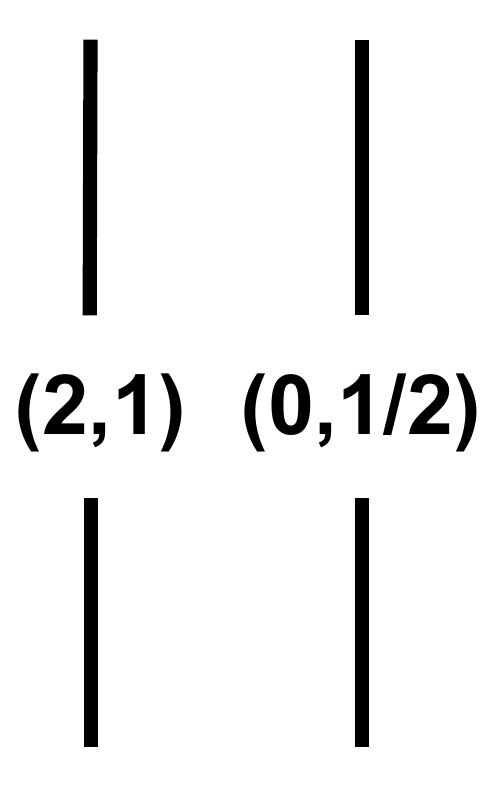}
\label{subfig:ex}
}
\\
\subfigure[]{
\includegraphics[scale=0.3]{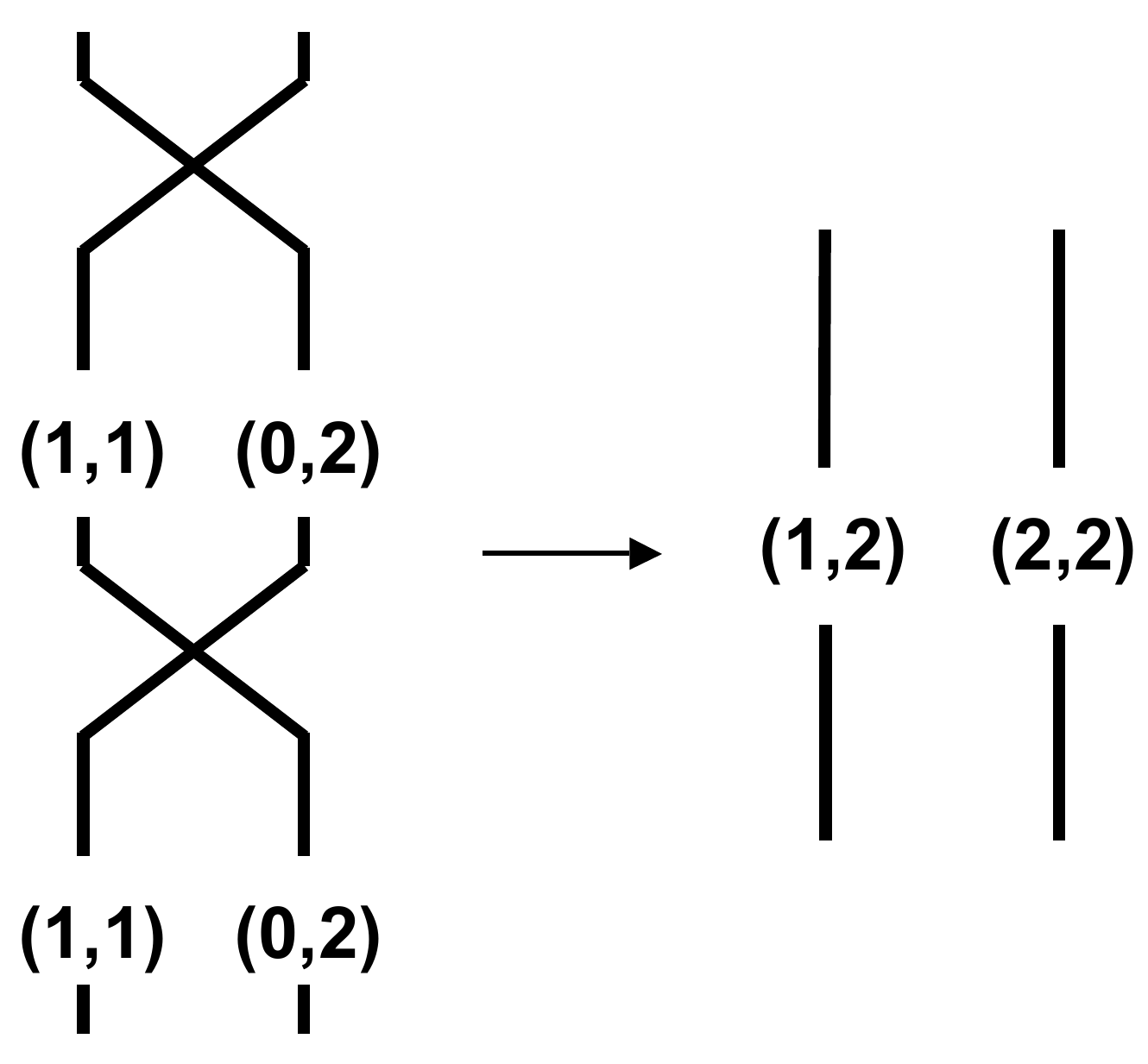}
\label{subfig:comp}
}
\caption{\subref{subfig:generalform} General form of the morphisms of $Hom(A^{\otimes 2},A^{\otimes 2})$ in $\mathcal{C}$; here, $g_1,g_2 \in Aff_+(\mathbb{R})$. \subref{subfig:ex} Two examples of morphisms \subref{subfig:comp} The composition of two morphisms (the order of composition is from top to bottom)}
\label{fig:morphismAn}
\end{figure}

The morphisms of $A^{\otimes n}$ can be represented graphically as string diagrams, as pictured in Figure \ref{fig:morphismAn} in the case $n=2$. The timelines are illustrated by the strings and the elements of $Aff_+(\mathbb{R})$ are the time-spans transformations on each timeline. The composition of two morphisms is also illustrated in this Figure, the composition order being from top to bottom of the diagram.

\begin{figure}
\centering
\includegraphics[scale=0.3]{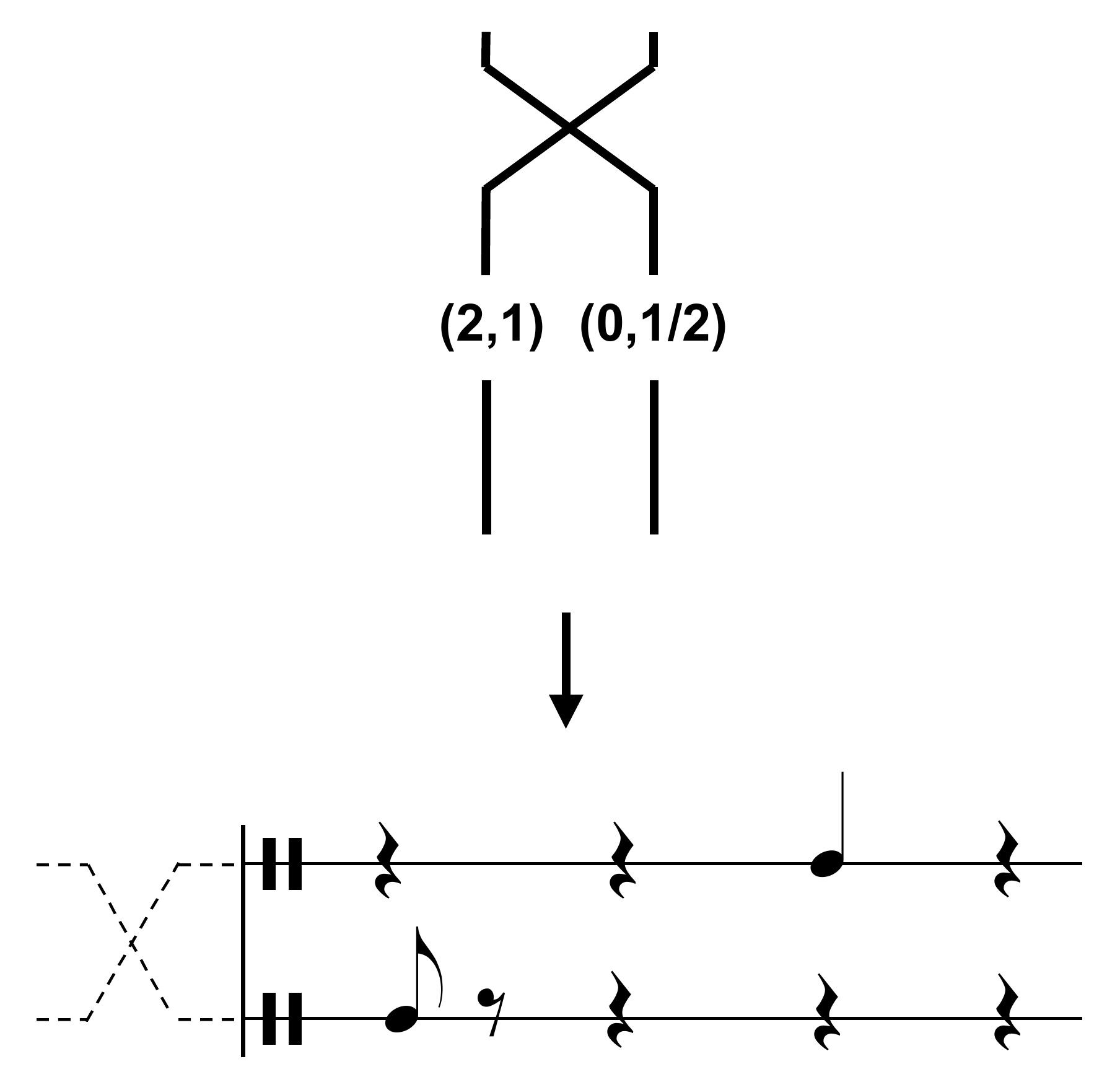}
\caption{The representation of a morphism of $A^{\otimes 2}$ as time-spans on two timelines. The additional information provided by the symmetry operation is represented by dashed lines between the two timelines.}
\label{fig:RepMorphism}
\end{figure}

The morphisms of $A^{\otimes n}$ are transformations of the time-spans, but, as elements of the image set of $A^{\otimes n}$ by a representable functor, they can be considered as the time-spans themselves. We can conveniently represent a morphism of the type $(g_1\otimes \cdots \otimes g_n) \circ \sigma$ (with $g_i  \in G$ and $\sigma$ is generated by the $\sigma_i$ morphisms), by time-spans $g_1, \cdots , g_n$ on timelines $1 \cdots n$, with the additional interexchange information provided by $\sigma$. An example of such a representation is pictured in the case $n=2$ in Figure \ref{fig:RepMorphism}.

\begin{figure}
\centering
\includegraphics[scale=0.5]{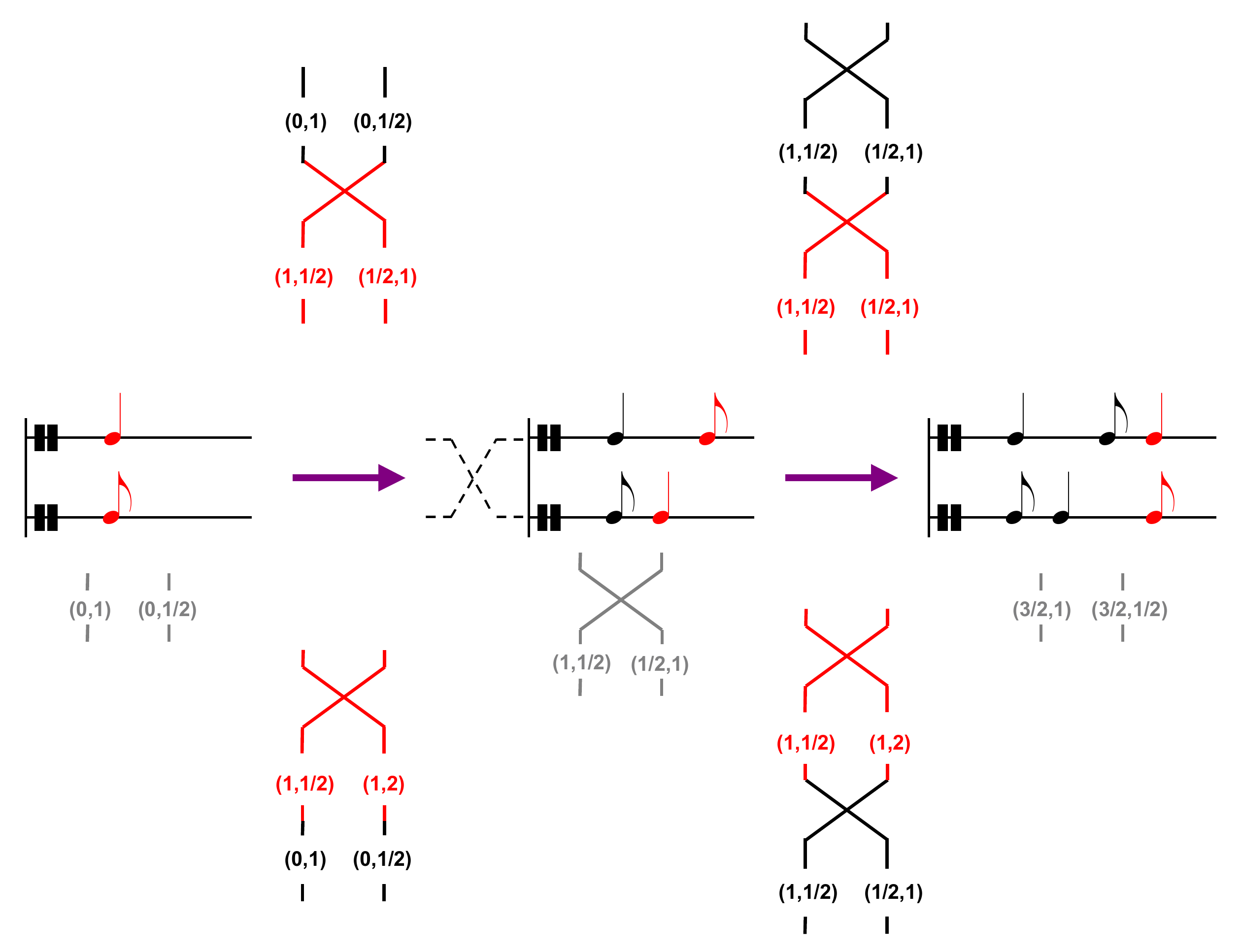}
\caption{Analysis of the rhythm of Figure \ref{fig:Rhythm}. The different time-spans result from the iterated covariant (top) or contravariant (bottom) composition of a single morphism, represented in red. The morphisms in gray are the morphisms corresponding to the highlighted red time-spans.}
\label{fig:RhythmAnalysis}
\end{figure}

With such a representation of morphisms, it is easy to analyze the rhythm of Figure \ref{fig:Rhythm}. The initial time-spans correspond to the morphism $((0,1) \otimes (0,\frac{1}{2}))$. If a covariant representable functor has been chosen, then the rhythm results from the iterated composition of this initial morphism by the morphism $((1,1) \otimes (\frac{1}{2},1)) \circ \gamma$. If on the other hand a contravariant representable functor has been chosen, then the rhythm results from the iterated composition of the initial morphism by the morphism $((1,\frac{1}{2}) \otimes (1,2)) \circ \gamma$. In both case, the presence of the symmetric operation $\gamma$ reflects the symmetry between the two timelines which is characteristic of this rhythm. The summary of this analysis is depicted graphically on Figure \ref{fig:RhythmAnalysis}.

\section{Encapsulated time-spans}

The construction introduced in the previous section allows the study of time-spans on multiple timelines. However, this model only deals with the transformation of a single time-span per timeline. This can be a limitation in some musical problems.

\begin{figure}
\centering
\includegraphics[scale=0.5]{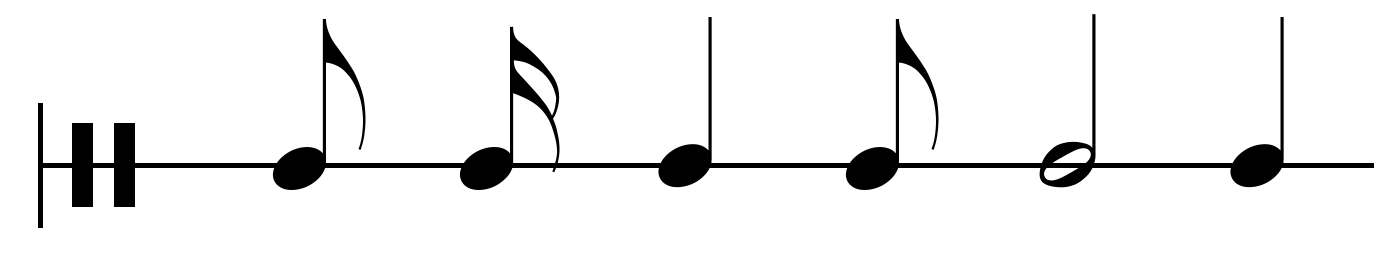}
\caption{A single timeline, dilating rhythm.}
\label{fig:RhythmDouble}
\end{figure}

Consider for example the rhythm presented in Figure \ref{fig:RhythmDouble}. The analysis of the successive time-spans using the group of transformations $Aff_+(\mathbb{R})$ would involve the successive left actions of the group elements $(\frac{1}{2},\frac{1}{2})$, $(-\frac{1}{4},4)$, $(\frac{11}{8},\frac{1}{2})$, $(-\frac{5}{4},4)$, etc. on the initial time-span $(0,\frac{1}{2})$. Obviously, this analysis misses completely the fact that the rhythm is generated by the successive dilation and translation of an initial unit cell composed the first quaver and semiquaver note. An analysis involving right actions would be more enlightening as it would consider the alternative right actions of the group elements $(1,\frac{1}{2})$ and $(1,4)$ on the initial time-span $(0,\frac{1}{2})$. However, we still need two different group elements to describe a single dilation of the initial cell.

Another example is provided by violin playing. Two notes (or more) can be player by a violin player in one bow. In that case, the two notes have their own time-spans but they are "encapsulated" in the time-span of the bow action itself. Note that if the player takes, for example, twice more time to perform the bow, the two notes will be dilated and translated correspondingly. In other words, a transformation applied to encapsulated time-spans will result in encapsulated, individually transformed, time-spans.

In mathematical terms, we would intuitively like to have a morphism $\Box$ taking a number $n$ of single time-spans from $n$ different timelines and "encapsulating" them into a "bracket" on a single timeline, or in other words a morphism $\Box : A^{\otimes n} \to A$. However, we can go further in this direction: indeed, we could also encapsulate both single time-spans and encapsulated time-spans, forming "brackets of brackets".

A solution for constructing such objects is the following. As before, we take the single category $\mathcal{G}$ with one object $A$ such that $Hom(A,A) \cong G$, where $G$ is a group. The construction is valid for any group $G$, though we will focus here on $G=Aff_+(\mathbb{R})$ with the multiplication operation given by

$$(t_1,\delta_1) \cdot (t_2,\delta_2) = (t_1+\delta_1 t_2, \delta_1 \delta_2).$$

We also assume that there exists a morphism $\Box : A^{\otimes 2} \to A$. This will give rise to brackets containing only two time-spans, but this can easily be generalized to any number of time-spans per bracket. We then construct the strict monoidal category $\mathcal{C}$ freely generated by $\mathcal{G}$ and by $\Box : A^{\otimes 2} \to A$. It can easily be seen that $Hom(A^{\otimes m},A^{\otimes n})=\emptyset$ for $m<n$ and that the morphisms of $Hom(A^{\otimes m},A^{\otimes n})$ for $m>n$ consist of forests of $n$ binary trees with $m$ total leaves, each branch of the trees (including the leaves) being decorated with an element of $G$. However we need to impose one relation on $\mathcal{C}$ to reflect the operations which can be applied to brackets of time-spans. Consider for example the morphism of $Hom(A^{\otimes 2},A)$ such as the one presented below, in which we have represented the morphism $\Box$ by the junction of two strings.

\vspace{0.2cm}
\begin{center}
\includegraphics[scale=0.3]{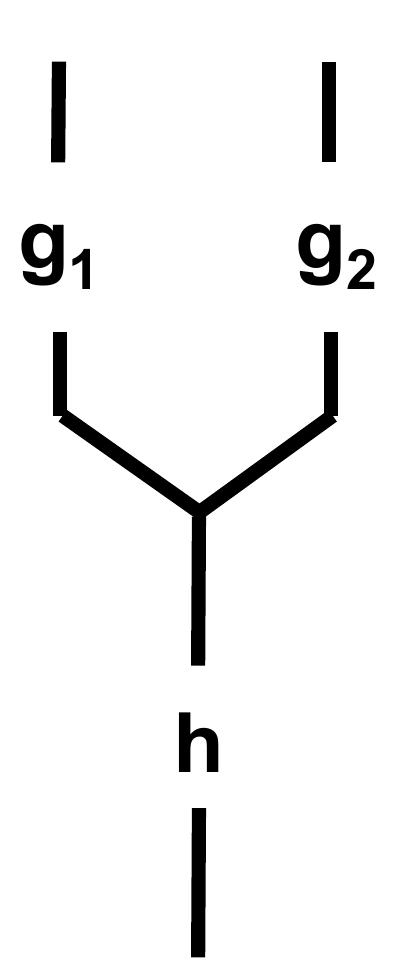}
\end{center}
\vspace{0.2cm}

The morphism $h$ applies a transformation of $G$ on the encapsulated time-spans $g_1$ and $g_2$. How should this morphism transform the time-spans ? Intuitively, a transformation of a bracket of time-spans should act as if the bracket itself was a time-span. For example, a global dilation of a bracket will induce a dilation of each of the encapsulated time-span, but the relative interval between the time-spans will stay the same. For the bracket $\Box \circ (g_1 \otimes g_2)$, the relative interval between the two time-spans is $g_1^{-1} \cdot g_2$. Therefore, if $g_1$ is transformed into $h_1$ and $g_2$ in $h_2$ we wish that $h_1^{-1} \cdot h_2 = g_1^{-1} \cdot g_2$. The straightforward way to do so is to impose the following relation :

$$ \mathcal{R} : h \circ \Box \circ (g_1 \otimes g_2) = \Box \circ (hÊ\cdot g_1 \otimes hÊ\cdot g_2)$$

Note that, contrary to the case where $G$ is used to describe transformations of single time-spans, we cannot use the right multiplication by $h$, for in that case the relative interval $g_1^{-1} \cdot g_2$ would change under transformation by $h$. In order to apply the equivalent action to a right multiplication on a bracket time-spans we will need to introduce new actions which we will describe later.

By quotienting $\mathcal{C}$ with the equivalence relation $\mathcal{R}$, we obtain the category $\mathcal{C}'$ in which $Hom(A^{\otimes m},A^{\otimes n})=\emptyset$ for $m<n$ and in which the morphisms of $Hom(A^{\otimes m},A^{\otimes n})$ for $m>n$ consist of forests of $n$ binary trees with $m$ total leaves, each leaf of the trees being decorated with an element of $G$.

\begin{figure}
\centering
\includegraphics[scale=0.4]{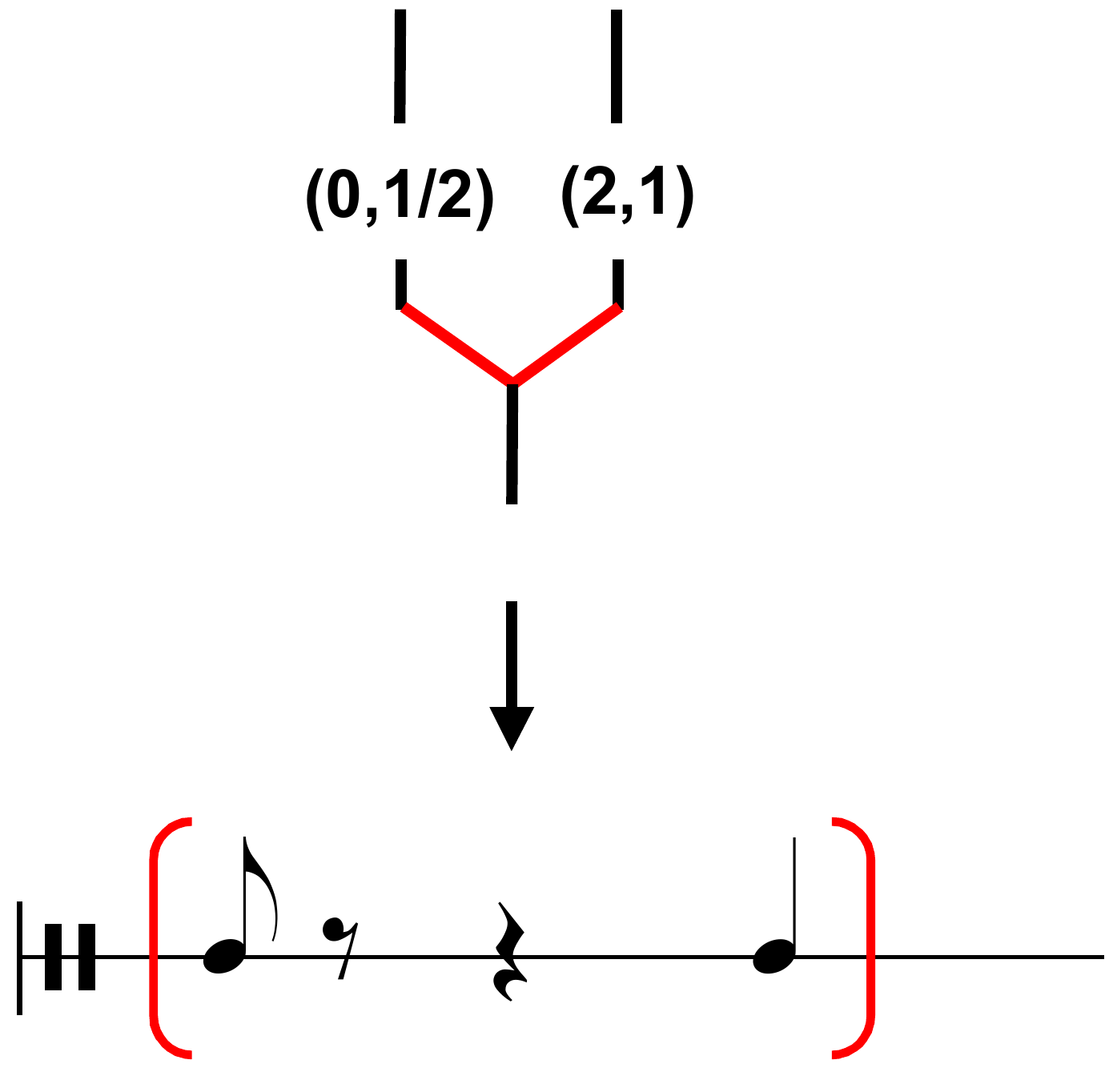}
\caption{The graphical representation of a morphism of $Hom(A^{\otimes 2},A)$}
\label{fig:BracketRepresentation}
\end{figure}

If we wish to consider the morphisms of $\mathcal{C}'$ as musical objects, we can, as previously, consider a representable covariant functor. For the reasons exposed above, we cannot choose this time
a contravariant functor, which would induce right multiplications on encapsulated time-spans. A $F=Hom(A^{\otimes m},-)$ functor will give rise to non-empty sets containing the morphisms of $Hom(A^{\otimes m},A^{\otimes n})$ (with $m>n$) which corresponds to $m$ total time-spans being encapsulated on $n$ timelines. In particular, each set $F(A^{\otimes n})$ is thus equipped with an action of $G^n$ which expresses transformations of time-spans and brackets for each time-line. A graphical representation of a morphism of $Hom(A^{\otimes m},A^{\otimes n})$ as a musical object can be obtained by drawing brackets around the corresponding notes on the staves. An example for a morphism of $Hom(A^{\otimes 2},A)$ is given in Figure \ref{fig:BracketRepresentation}.

\begin{figure}
\centering
\includegraphics[scale=0.6]{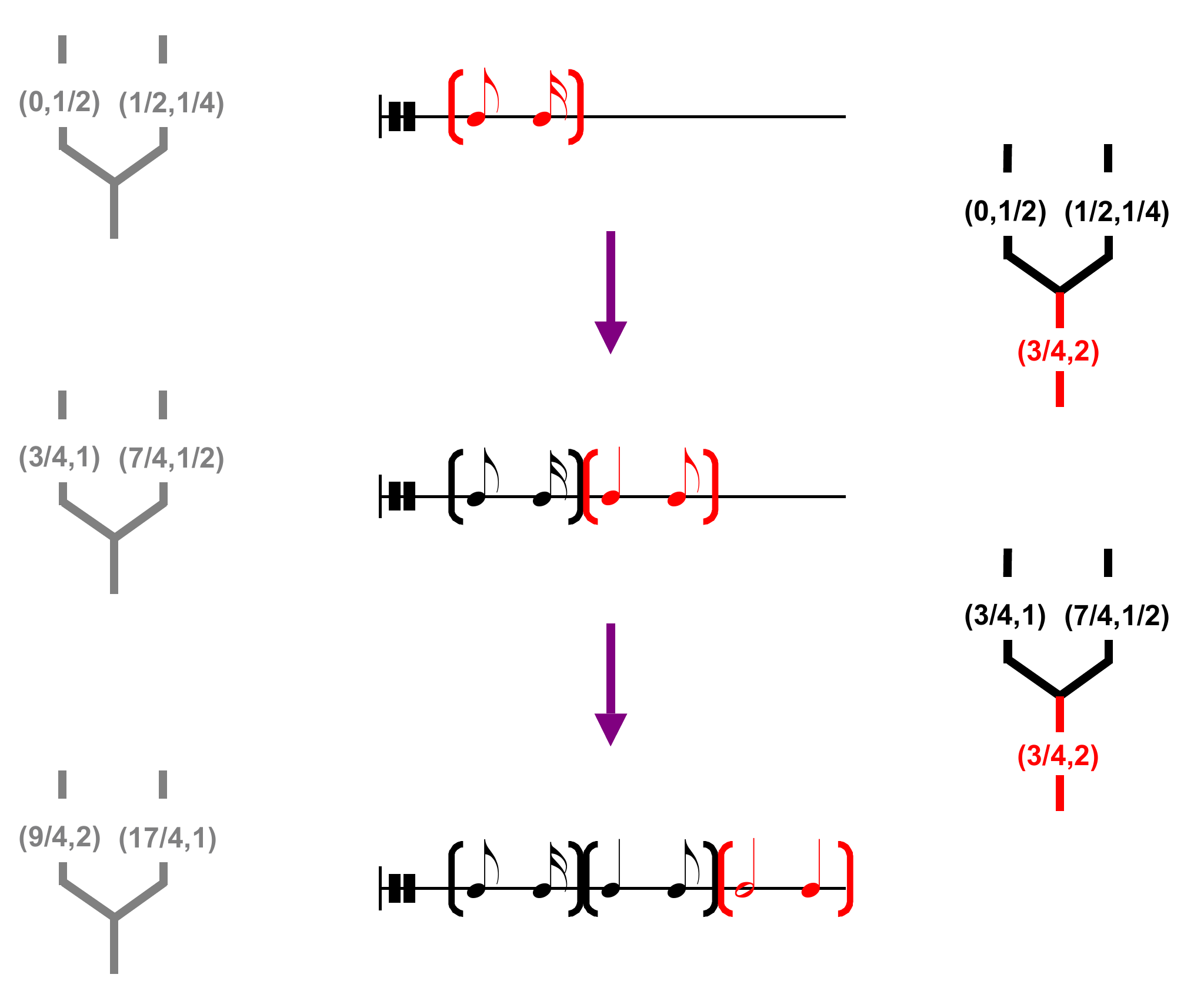}
\caption{Analysis of the rhythm of Figure \ref{fig:RhythmDouble}. The different time-spans result from the iterated covariant composition of a single morphism, represented in red. The morphisms in gray are the morphisms corresponding to the highlighted red encapsulated time-spans.}
\label{fig:BracketAnalysis}
\end{figure}

We are now in position to analyze the rhythm of Figure \ref{fig:RhythmDouble}. The initial cell consisting of the first quaver and semiquaver notes corresponds to the morphism $\Box \circ ( (0,\frac{1}{2}) \otimes (\frac{1}{2},\frac{1}{4}))$. The iterated application of $(\frac{3}{4},2)$ on this morphism will give rise to the successive time-spans composing the rhythm, as represented graphically on Figure \ref{fig:BracketAnalysis}. Thus we obtain a clear analysis of this rhythm using only one transformation, which we couldn't do by just considering the group $G=Aff_+(\mathbb{R})$.

A question remains: how can we define actions on brackets analogous to right actions for single time-spans ? We need here to reexamine the notion of right action which we described in the previous section. We have seen that given a group-as-category $G$ and a covariant representable functor, we obtain a set $S$, a canonical bijection $\chi : G \to S$ and a canonical left group action of an element $g$ of $G$ on an element $s$ of $S$ as

$$g \cdot s = \chi(g \cdot \chi^{-1}(s))$$
Notice that we can choose the element of $G$ to be of the form $u \cdot g \cdot u^{-1}$, $u \in G$ and we have

$$g \cdot s = \chi(u \cdot g \cdot u^{-1} \cdot \chi^{-1}(s))$$

The interpretation of this action is, by analogy with well-known cases in linear algebra for matrices, an action of the group element $g$ after a "change of basis" induced by the group element $u$.
In particular, we can choose to take $u=\chi^{-1}(s)$. We then obtain the action

$$g \cdot s = \chi(\chi^{-1}(s) \cdot g \cdot (\chi^{-1}(s))^{-1} \cdot \chi^{-1}(s)) = \chi(\chi^{-1}(s) \cdot g)$$
which is the right action corresponding to right multiplication. The particular interpretation we give in this case is that, by the "change of basis" $u=\chi^{-1}(s)$, we have made $s$ the new reference to which a transformation $g$ can be applied. It is worthwhile to note that this exactly what Lewin advocated when he devised his non-commutative GIS for time-spans: in \cite{lewin}, he says \textit{"(...) any time-span has the potential for becoming locally referential, or behaving as if it were"}.

Note that this point of view can also be applied to chords. Indeed the famous contextual $P$, $L$ and $R$ operations acting on major and minor triads can be understood as the classical non-contextual inversion operators applied in the referential frame of the chord, as exposed for example in \cite{fiore3}.

To transpose this point of view to brackets of time-spans we will take an example familiar to users of graphic softwares. Suppose we have two squares A and B, which we wish to scale. We can apply a global dilation, as pictured below, with respect to a given referential frame. Notice that the squares B' and A' are in the same proportions as are the squares A and B. 

\vspace{0.2cm}
\begin{center}
\includegraphics[scale=0.4]{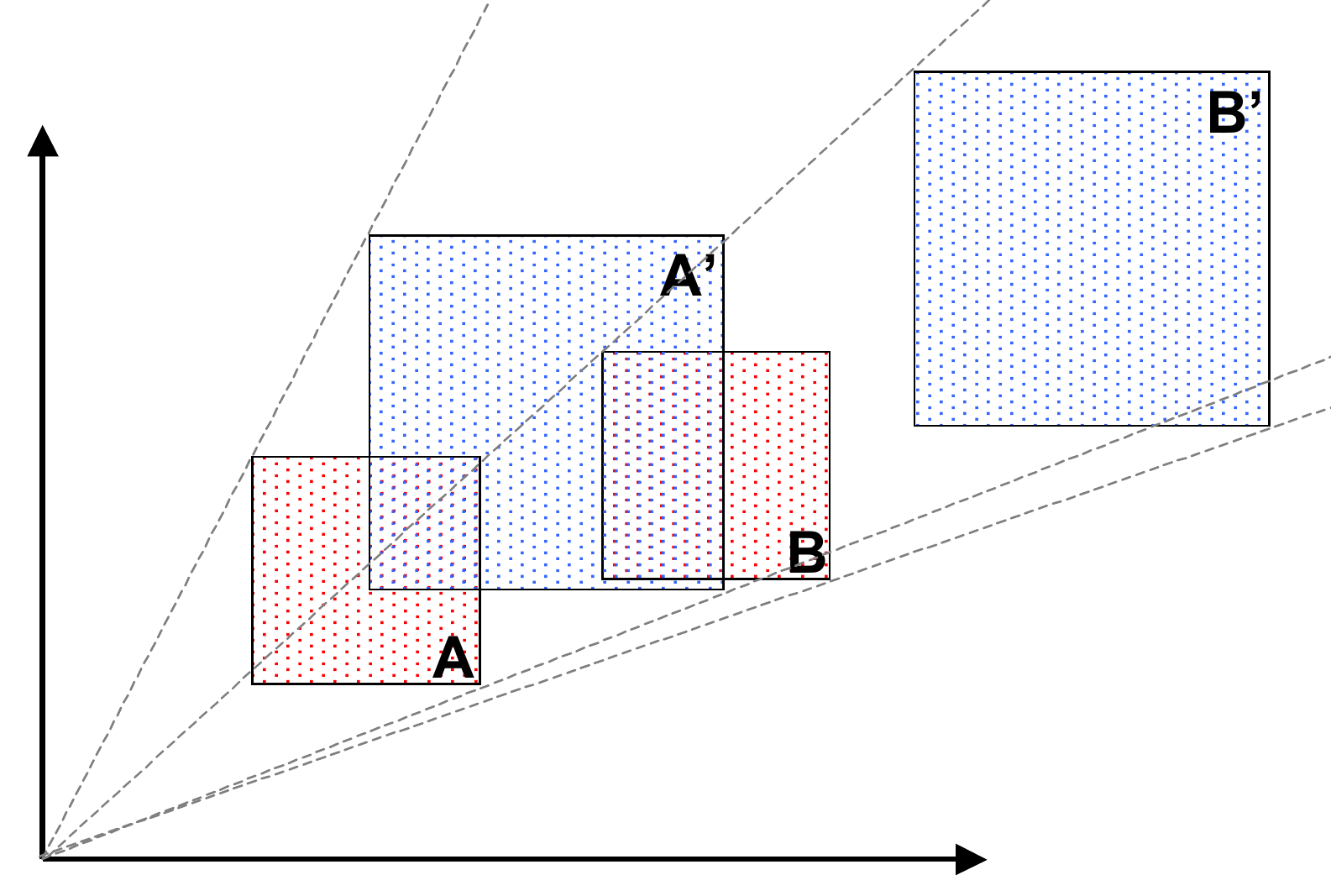}
\end{center}
\vspace{0.2cm}

If we dilate each square individually, we lose the relative proportion between squares, as pictured below.

\begin{center}
\includegraphics[scale=0.4]{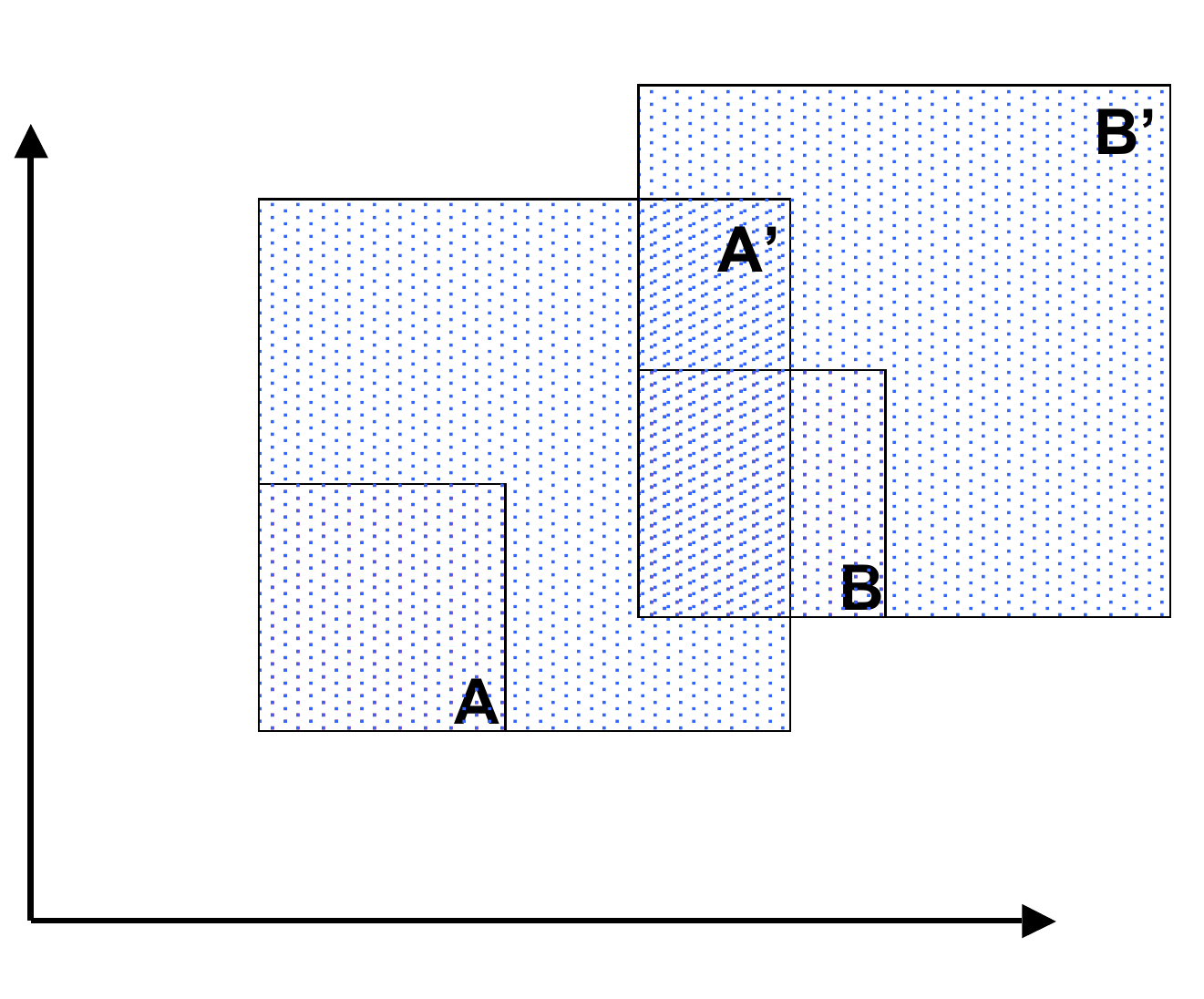}
\end{center}
\vspace{0.2cm}

Instead, we can translate the referential frame to the corner of square A, dilate both squares with a global transformation in this new referential, then move it back to its original position, as pictured below (this is quite often the transformation applied implicitly when objects are "grouped" in graphic softwares).

\begin{center}
\includegraphics[scale=0.4]{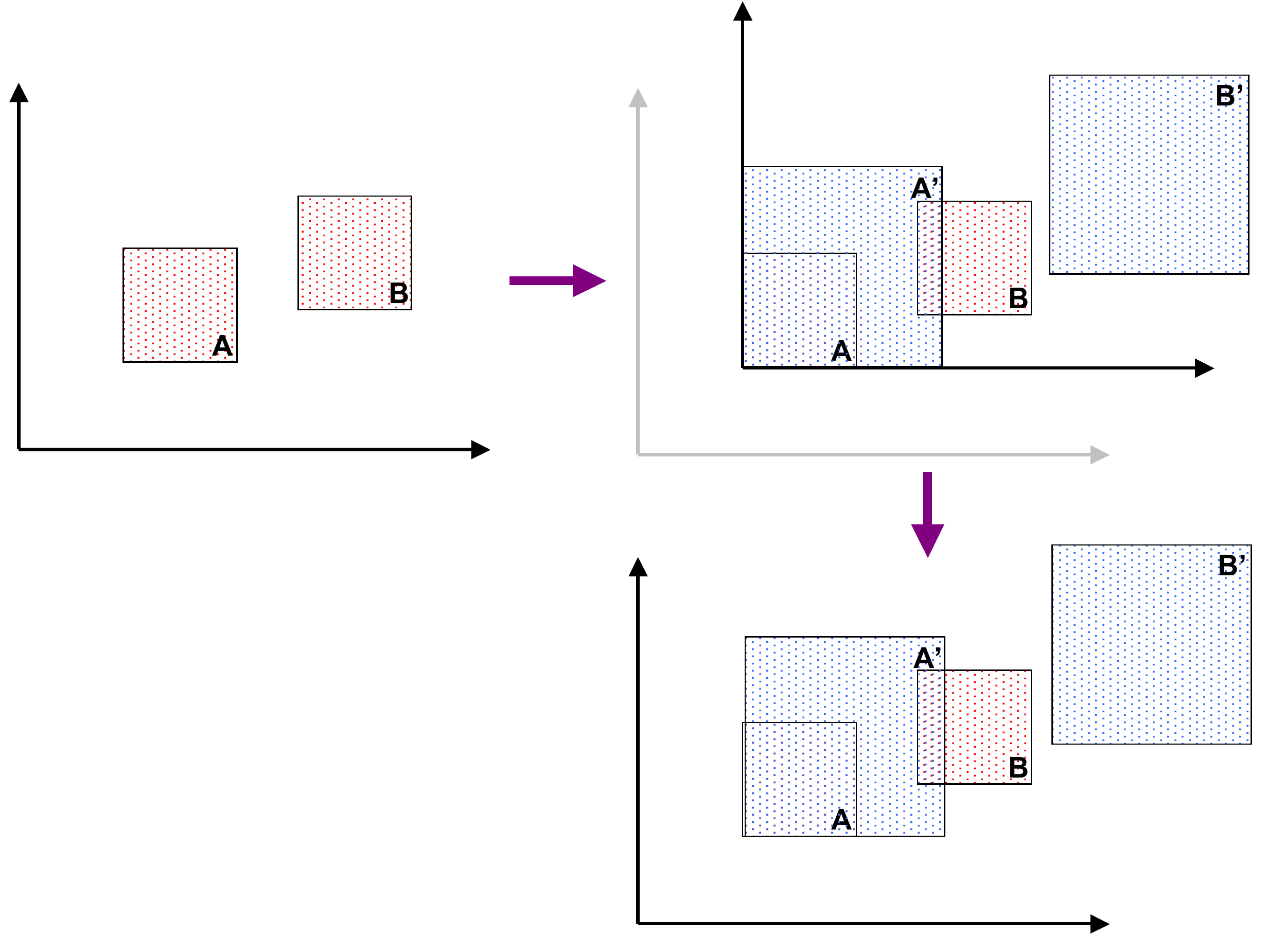}
\end{center}
\vspace{0.2cm}

In this transformation, the square A has become the new reference in which the dilation is applied, in a analogous manner to the case of single time-spans described above. Of course, we could also have chosen the square B as the new reference.

The application for brackets of time-spans is similar to this graphic example. Assume for example that we have a bracket of two time-spans which correspond to the morphism $\mu=\Box \circ (g_1 \otimes g_2)$, and that we want to apply a transformation $h$, \textit{relative to the first time-span}. We will then compose $\mu$ with $g_1 h g_1^{-1}$ to obtain the bracket of time-spans corresponding to the morphism $\mu'=\Box \circ ((g_1 h) \otimes (g_1 h g_1^{-1} g_2))$. Notice that the relative interval $g_1^{-1}g_2$ has been preserved in this transformation. If we wish to apply another transformation $k$ to this new bracket, relative to the first time-span, we would then have to compose the morphism $\mu'$ with 
$(g_1h) k (g_1h)^{-1}$. In this way, we can also analyze the rhythm of Figure \ref{fig:RhythmDouble} by the iterated composition of the initial morphism $\Box \circ ( (0,\frac{1}{2}) \otimes (\frac{1}{2},\frac{1}{4}))$ with the morphism $(\frac{3}{2},2)$, \textit{relative to the first time-span}. Notice that $\frac{3}{2}$ is precisely the total length of the initial cell, as measured with the first time-span as reference, while the value $2$ reflects the dilation which is applied. Hence this contextual transformation precisely expresses
the translation by one unit-cell and dilation which is observed in this rhythm. While our example refers to a bracket of two time-spans, this kind of transformation can be applied to any bracket, by selecting the group element of one of the leaves in the corresponding morphism as the reference and applying the corresponding morphism.

Finally, one can note that we have not introduced any braiding to express the possible interexchange of brackets belonging to different timelines. This can be done by constructing the strict \textit{symmetric} monoidal category $\mathcal{C}$ freely generated by $\mathcal{G}$ and by $\Box : A^{\otimes 2} \to A$. As before we will need to quotient this category by the relation
$$ \mathcal{R} _1: h \circ \Box \circ (g_1 \otimes g_2) = \Box \circ (hÊ\cdot g_1 \otimes hÊ\cdot g_2)$$

However we also need to introduce new relations to express the result of echanging a bracket and a time-span or two different brackets, namely 

$$ \mathcal{R}_2 : \gamma \circ (id \otimes \Box) = (\Box \otimes id) \circ (id \otimes \gamma) \circ (\gamma \otimes id)$$
$$ \mathcal{R}_3 : \gamma \circ (\Box \otimes id) = (id \otimes \Box) \circ (\gamma \otimes id) \circ (id \otimes \gamma)$$
$$ \mathcal{R}_4 : (\gamma \otimes \gamma) \circ (\Box \otimes \Box) = (\Box \otimes \Box) \circ (id \otimes \gamma \otimes id) \circ (\gamma \otimes \gamma) \circ (id \otimes \gamma \otimes id)$$

Quotienting $\mathcal{C}$ by the relations $\mathcal{R}_{1\cdots4}$ will generate a category $\mathcal{C}'$ in which the reader can easily verify that the morphisms of $Hom(A^{\otimes m},A^{\otimes n})$ for $m>n$ can be put in the form of a successive composition of a permutation $\sigma$ of the $m$ leaves, followed by group elements of $G$, followed by bracketing operations. The use of these morphisms as musical objects and their corresponding transformations is straightforward by following the same approach as above.

\section{Conclusion}

In this paper we have generalized the known transformational theory of time-spans to both multiple time-lines and brackets of time-spans using the abstract framework of monoidal categories. In particular we have shown that such a framework provides useful tools for the analysis of complex rhythms which would be otherwise untractable using only the group $G=Aff_+(\mathbb{R})$ introduced by Lewin.

Musical transformational theory has been mainly applied to chords, for example major and minor chords, but seventh chords as well. It would therefore be interesting to see if monoidal categories arise in such situations and what would be the benefit of such tools for musical analysis.


\label{lastpage}

\end{document}